# Character Average of Second and Fourth Powers of Dirichlet L-series at Unity

VIVEK V. RANE .
Department of Mathematics ,
The Institute of Science,
15, Madame Cama Road,
Mumbai-400 032,
India .
v_v_rane@yahoo.co.in

**Abstract** : For a Dirichlet character modulo an integer $q \geq 3$, we use a highly simple elementary method to give an asymptotic formula for $\sum_{\chi \neq \chi_0 (\bmod q)} |L(1,\chi)|^4$, where $\chi_0 (\bmod q)$ is the principal character . This result seems to be new . We also obtain an asymptotic formula for $\sum_{\chi \neq \chi_0 (\bmod q)} |L(1,\chi)|^2$ , using power series expansion of $\zeta(s, 1+\alpha)$ at $s = 1$.

# Character Average of Second and Fourth Powers of Dirichlet L-series at Unity

VIVEK V. RANE .

Department of Mathematics ,
The Institute of Science,
15, Madame Cama Road,
Mumbai-400 032,
India .
v_v_rane@yahoo.co.in

**Introduction** : Let $s=\sigma+it$ be a complex variable , where $\sigma,t$ are real . For an integer $q\geq 1$ , let $\chi(\bmod q)$ denote a Dirichlet character and let $\chi_0(\bmod q)$ denote the principal character . Let $L(s,\chi)$ be the corresponding Dirichlet L-series so that $L(s,\chi)=\sum_{n\geq 1}\chi(n)n^{-s}$ for Re $s>1$ ; and its analytic continuation . For $0<\alpha\leq 1$ , let $\zeta(s,\alpha)$ be the Hurwitz zeta function defined by

$\zeta(s,\alpha)=\sum_{n\geq 0}(n+\alpha)^{-s}$ for Re $s>1$ ; and its analytic continuation .

For $0\leq\alpha\leq 1$ , let $\zeta_1(s,\alpha)=\sum_{n\geq 1}(n+\alpha)^{-s}$ for Re $s>1$ ; and its analytic continuation . Let $\zeta(s)$ denote the Riemann zeta function . Thus

$\zeta(s,1)=\zeta_1(s,0)=\zeta(s)$ . Also note that $L(s,\chi)=q^{-s}\sum_{a=1}^{q}\chi(a)\zeta(s,\frac{a}{q})$ .

For an integer $n\geq 1$ , we shall write $a(n)=a_q(n)=\sum_{uv=n}1$ with $1\leq u$ , $v\leq q$ , for the given integer $q\geq 3$ . Note that $a(n)\leq d(n)$ , where $d(n)$ is the divisor function .

In what follows , $\phi(n)$ shall denote Euler's totient function and $\mu(n)$ shall stand for Moebius function and p shall denote a prime number. Also $\gamma$ will mean Euler's constant . For $u$ real , $[u]$ shall stand for its integral part . In what follows , $\sum_a{}^{'} \left( or \sum_{a=1}^{q'} \right)$ shall mean the summation $\sum_{a=1}^{q}$ with the restriction $(a,q)=1$ . The summation $\sum_b{}^{'}$ , $\sum_\ell{}^{'}$ and $\sum_k{}^{'}$ will have similar meanings . In what follows , the 0 and << constants will be absolute , unless stated otherwise .

Effectively , we shall be using only three facts namely ,

1) Euler's summation Formula

2) the boundedly convergent series expression on the unit interval ,

$u-\frac{1}{2}=-\sum_{|n|\geq 1}\frac{e^{2\pi inu}}{2\pi in}$ with equality for non-integral u .

3) Power series (in $\alpha$ ) of $\zeta_1(s,\alpha)$ .

Our object is to prove the following theorems .

**Theorem** 1 : For an integer $q \geq 3$ , we have

$$\sum_{\chi\neq\chi_0 (\bmod q)} |L(1,\chi)|^4 = \phi(q)\cdot \sum_{\substack{n=1\\(n,q)=1}}^{q} \frac{d^2(n)}{n^2} + 0\left(\frac{\phi(q)}{\sqrt{q}}\log^4(q+2)\right) ,$$

where $\chi_0 (\bmod q)$ is the principal character and $d(n)$ is the number of divisors of $n$ .

**Remark** : Our Theorem above for $s=1$ , seems to be the best result in this direction . Our method is highly simple and elementary .

Before stating our Theorem 2 , we shall need a Lemma .

**Lemma** : Let $\zeta_1(s,\alpha)=\zeta(s,\alpha)-\alpha^{-s}$ . Then for any complex number $s$ ,

we have $\zeta_1(s,\alpha)=\zeta(s)+\sum_{n\geq 1} a_n(s)\alpha^n$ , where $a_n(s)=\frac{(-1)^n}{n!}s(s+1)........(s+n-1)\zeta(s+n)$ in the open disc $|\alpha|<1$ of $\alpha$ - plane and also at $\alpha=1$.

**Corollary** : We have

$$\lim_{s\to 1}\left(\zeta_1(s,\alpha)-\tfrac{1}{s-1}\right)=\lim_{s\to 1}(\zeta(s)-\tfrac{1}{s-1})+\lim_{s\to 1}\sum_{n\geq 1}\tfrac{(-1)^n}{n!}s(s+1)......(s+n-1)\zeta(s+n)\alpha^n=\gamma+\sum_{n\geq 1}(-1)^n\zeta(n+1)\alpha^n$$

**Remark** : Our lemma above is Theorem 1 of author [2] . We shall write

$\lim_{s\to 1}\left(\zeta_1(s,\alpha)-\tfrac{1}{s-1}\right)=\sum_{n\geq 0} a_n\alpha^n$ , where $a_0=\gamma$ and $a_n=(-1)^n\zeta(n+1)$ for $n\geq 1$.

Let $\left(\lim_{s\to 1}(\zeta_1(s,\alpha)-\tfrac{1}{s-1})\right)^2=\sum_{n\geq 0} c_n\alpha^n$ , where $c_n=\sum_{i+j=n} a_i a_j$ for $n\geq 0$.

Incidentally , it is to be noted that for $s>0$ and for $\alpha$ in the unit interval [0,1] and for any integer $N\geq 1$, we have finite Taylor series namely ,

$$\zeta_1(s,\alpha)=\sum_{n=0}^{N-1} a_n(s)\alpha^n+\tfrac{(-1)^N}{N!}s(s+1)..........(s+N-1)\cdot\zeta_1(s+N,\theta_1\alpha)\alpha^N$$

for some $\theta_1=\theta_1(s,N)$ where $0<\theta_1<1$ so that for $s>0$,

$\zeta_1(s,\alpha)=\sum_{n=0}^{N-1} a_n(s)\alpha^n+0\left(\frac{s(s+1)..........(s+N-1)}{N!}\cdot\zeta(s+N)\alpha^N\right)$ , where the O-constant is 1 .

In particular , letting $s\to 1$, we have $\lim_{s\to 1}\left(\zeta_1(s,\alpha)-\tfrac{1}{s-1}\right)=\sum_{n=0}^{N-1} a_n\alpha^n+0\left(\zeta(N+1)\alpha^N\right)$ ,

where the O-constant is 1 . From this , we get the finite Taylor series of

$g(\alpha)=\left(\lim_{s\to 1}\left(\zeta_1(s,\alpha)-\tfrac{1}{s-1}\right)\right)^2$ as follows : For any integer $N\geq 1$,

$\left(\lim_{s\to 1}(\zeta_1(s,\alpha)-\tfrac{1}{s-1})\right)^2=\sum_{n=0}^{N-1} c_n\alpha^n+\tfrac{\alpha^N}{N!}g^{(N)}(\theta_2\alpha)$ for some $\theta_2=\theta_2(N)$ with

$0<\theta_2<1$, where $g^{(N)}(\alpha)=\tfrac{d^N}{d\alpha^N}g(\alpha)$.

Thus $g(\alpha)=\sum_{n=0}^{N-1} c_n\alpha^n+0(\alpha^N)$, where the O-constant is dependent on N .

Next, we state our Theorem 2 .

**Theorem 2** : We have for $q \geq 3$ ,

$$\sum_{\chi \neq \chi_0 (\bmod q)} |L(1,\chi)|^2 = \phi(q) \sum_{\substack{a=1 \\ (a,q=1)}}^{q} \frac{1}{a^2} - \frac{\phi^2(q)}{q^2} \left( \log q + \sum_{p/q} \frac{\log p}{p-1} \right)^2 + \frac{\gamma^2 \phi^2(q)}{q^2} + 0\left( \frac{d(q)\phi(q)}{q^2} \right)$$

$$+ \frac{2\phi(q)}{q^2} \cdot \sum_{n \geq 1} \frac{(-1)^n \zeta(n+1) q^{1-n}}{n} \cdot \sum_{k/q} \mu(\tfrac{q}{k}) \left( B_n(k+1) - B_n(1) \right)$$

$$+ \frac{\phi(q)}{q^2} \cdot \sum_{n \geq 0} \frac{c_n}{n+1} \cdot q^{-n} \cdot \sum_{k/q} \mu(\tfrac{q}{k}) \left( B_{n+1}(k+1) - B_{n+1}(1) \right)$$

, where $B_n(x)$ is Bernoulli polynomial of degree $n$; and for $n \geq 0$, $c_n = \sum_{i+j=n} a_i a_j$ , with $a_0 = \gamma$ and $a_n = (-1)^n \zeta(n+1)$ when $n \geq 1$.

Remark : In the light of finite Taylor series (with remainders) for $\lim_{s \to 1} \left( \zeta_1(s,\alpha) - \frac{1}{s-1} \right)$ and $\left( \lim_{s \to 1} (\zeta_1(s,\alpha) - \frac{1}{s-1}) \right)^2$, the last two terms involving infinite series in the statement of Theorem 2 above can be replaced by asymptotic expansions with error term $0(q^{-N})$ for any positive integer N. Also is the case for the term $0\left( \frac{d(q)\phi(q)}{q^2} \right)$ .

There are many papers on the estimation of $\sum_{\chi \neq \chi_0 (\bmod q)} |L(1,\chi)|^2$. However , the latest result in the case of $\sum_{\chi \neq \chi_0 (\bmod q)} |L(1,\chi)|^2$ has been given by S. Kanemitsu , Y. Tanigawa , M. Yoshimoto , W. Zhang [1] , which states that

$$\sum_{\chi \neq \chi_0 (\bmod q)} |L(1,\chi)|^2 = \zeta(2)\phi(q) \prod_{p/q} (1 - p^{-2}) - \frac{\phi^2(q)}{q^2} \left( \log q + \sum_{p/q} \frac{\log p}{p-1} \right)^2$$

$$+ \frac{\phi^2(q)}{q^2} \left( \gamma^2 - 2\gamma_1 - 2\zeta(2) \right) + \frac{2\phi(q)}{q^2} R(q) \ ,$$

where $R(q) = \sum_{d/q} \mu(\tfrac{q}{d}) \int_0^1 \overline{B_1}(du) \sum_{n \geq 1} \frac{H_n}{(n+u)^2} du$ with $\overline{B_1}(u) = B_1(u - [u])$,

$B_1(u)$ being Bernoulli polynomial of degree 1 ; $H_n = \sum_{k=1}^{n} \frac{1}{k}$

and $\gamma_1$ is the constant defined by $\zeta(s) = \frac{1}{s-1} + \gamma + \gamma_1(s-1) + \ldots\ldots\ldots$

**<u>Proof</u> <u>of</u> <u>Theorem</u> 1** : For a Dirichlet character $\chi \neq \chi_0 (\bmod q)$ and for Re $s > 1$ ,

we have $L^2(s,\chi) = \sum_{\ell}' \chi(\ell) \sum_{r \equiv \ell (\bmod q)} \frac{d(r)}{r^s}$ .

We shall write $Z_2(s,\ell,q) = \sum_{\substack{r \geq 1 \\ r \equiv \ell (\bmod q)}} \frac{d(r)}{r^s}$ for Re $s > 1$ ; and its analytic continuation .

Henceforth in what follows , in any congruence relation , for brevity we shall write $'q'$ in place of $'\bmod q'$ .

Let $\ell$ be an integer with $1 \leq \ell \leq q$ and $(\ell, q) = 1$ .

Thus for $\sigma > 1, Z_2(s,\ell,q) = \sum_{\substack{r \geq 1 \\ r \equiv \ell(q)}} \frac{d(r)}{r^s} = \sum\sum_{mn \equiv \ell(q)} (mn)^{-s}$

$$= \sum_{ab}' \sum_{\equiv \ell(q)}' \left( \sum_{m \equiv a(q)} m^{-s} \right) \left( \sum_{n \equiv b(q)} n^{-s} \right) = \sum_{ab}' \sum_{\equiv \ell(q)}' q^{-s} \zeta(s, \tfrac{a}{q}) q^{-s} \zeta(s, \tfrac{b}{q}) \ .$$

Thus for $(\ell, q) = 1$, and for $\chi \neq \chi_0 (\bmod q)$ and for $s \neq 1$,

we have $L^2(s,\chi) = \sum_{\ell}' \chi(\ell) \sum_{ab}' \sum_{\equiv \ell(q)}' q^{-s} \zeta(s, \frac{a}{q}) q^{-s} \zeta(s, \frac{b}{q})$ .

We shall show that for $s \neq 1$ and for $\chi \neq \chi_0 (\bmod q)$ ,

$$L^2(s,\chi) = \sum_{\ell}' \chi(\ell) \sum\sum_{ab \equiv \ell(q)}' q^{-s} \zeta(s, \tfrac{a}{q}) \cdot q^{-s} \zeta(s, \tfrac{b}{q})$$
$$= \sum_{\ell}' \chi(\ell) \sum\sum_{ab \equiv \ell (q)}' q^{-s} \left( \zeta(s, \tfrac{a}{q}) - \tfrac{1}{s-1} \right) \cdot q^{-s} \left( \zeta(s, \tfrac{b}{q}) - \tfrac{1}{s-1} \right).$$

From the expression $\zeta(s,\alpha) = \sum_{n \geq 0} (n+\alpha)^{-s}$ for $\sigma > 1$ , on using Euler's summation formula , for arbitrary $x > 0$ and for $\sigma > 1$ and for $0 < \alpha \leq 1$ , we have

$$\zeta(s,\alpha)=\sum_{0\le n\le x-\alpha}(n+\alpha)^{-s}+\tfrac{x^{1-s}}{s-1}+(x-\alpha-[x-\alpha]-\tfrac{1}{2})x^{-s}-s\int_x^\infty\tfrac{(u-\alpha-[u-\alpha]-\frac{1}{2})}{u^{s+1}}du$$

Note that $x-\alpha>-1$ . Here empty sum is treated as zero .

This expression is valid for $\sigma>0$ and for $s\neq 1$ . We choose $x=1$ . This gives for Re $s>0$ and $s\neq 1$,

$$\zeta(s,\alpha)=\alpha^{-s}+\tfrac{1}{s-1}+(\tfrac{1}{2}-\alpha)-s\int_1^\infty\tfrac{(u-\alpha-[u-\alpha]-\frac{1}{2})}{u^{s+1}}du\ .$$

(Incidentally , this gives $\lim_{s\to 1}\left(\zeta(s,\alpha)-\tfrac{1}{s-1}\right)=\alpha^{-1}+\left(\tfrac{1}{2}-\alpha\right)-\int_1^\infty\tfrac{\left(u-\alpha-[u-\alpha]-\frac{1}{2}\right)}{u^2}du$ ) .

In what follows , s is real positive unless stated otherwise .

For $s>0$ , $s\int_1^\infty\tfrac{(u-\alpha-[u-\alpha]-\frac{1}{2})}{u^{s+1}}du<<s\int_1^\infty\tfrac{du}{u^{s+1}}<<1$ , where << - constant is absolute .

Thus for $0<s\neq 1$ and for $0<\alpha\le 1$ , we have

$$\zeta(s,\alpha)=\alpha^{-s}+\tfrac{1}{s-1}+(\tfrac{1}{2}-\alpha)+0(1)$$

$=\alpha^{-s}+\tfrac{1}{s-1}+\phi(s,\alpha)$ , where $\phi(s,\alpha)=0(1)$ .

Thus $\zeta(s,\tfrac{a}{q})\zeta(s,\tfrac{b}{q})=(\tfrac{q^s}{a^s}+\tfrac{1}{s-1}+\phi(s,\tfrac{a}{q}))(\tfrac{q^s}{b^s}+\tfrac{1}{s-1}+\phi(s,\tfrac{b}{q}))$

$=((\tfrac{q^s}{a^s}+\phi(s,\tfrac{a}{q}))(\tfrac{q^s}{b^s}+\phi(s,\tfrac{b}{q}))+\tfrac{1}{s-1}(\tfrac{q^s}{a^s}+\phi(s,\tfrac{a}{q}))+\tfrac{1}{s-1}(\tfrac{q^s}{b^s}+\phi(s,\tfrac{b}{q}))+\tfrac{1}{(s-1)^2}$

Consider $\sum_{ab\equiv\ell(q)}'\sum'\tfrac{1}{s-1}(\tfrac{q^s}{a^s}+\phi(s,\tfrac{a}{q}))=\tfrac{1}{s-1}\sum_{ab\equiv\ell(q)}'\sum'\tfrac{q^s}{a^s}+\tfrac{1}{s-1}\sum_{ab\equiv\ell(q)}'\sum'\phi(s,\tfrac{a}{q})$

$=\tfrac{1}{s-1}\sum_a'\tfrac{q^s}{a^s}\sum_{\substack{b\\ ab\equiv\ell(q)}}'1+\tfrac{1}{s-1}\sum_{a=1}^{q}{}'\phi(s,\tfrac{a}{q})\sum_{\substack{b=1\\ ab\equiv\ell(q)}}^{q}{}'1=\tfrac{1}{s-1}\sum_a'\tfrac{q^s}{a^s}+\tfrac{1}{s-1}\sum_a'\phi(s,\tfrac{a}{q})$ ,

which is independent of $\ell$ .

Hence for $s \neq 1$ and for $\chi \neq \chi_0 (\operatorname{mod} q)$, we have $\sum_{\ell}' \chi(\ell)(\sum_{ab}' \sum_{\equiv \ell(q)}' \frac{1}{s-1}(\frac{q^s}{a^s} + \phi(s, \frac{a}{q})) = 0$ .

Similarly , $\sum_{\ell}' \chi(\ell) \left( \sum_{ab\equiv}' \sum_{\ell(q)}' \frac{1}{s-1}(\frac{q^s}{b^s} + \phi(s, \frac{b}{q})) \right) = 0$ for $s \neq 1$ and for $\chi \neq \chi_0 (\operatorname{mod} q)$ .

Similarly , we have $\sum_{\ell}' \chi(\ell) \sum_{ab\equiv}' \sum_{\ell(q)}' \frac{1}{(s-1)^2} = 0$ for $s \neq 1$ and for $\chi \neq \chi_0 (\operatorname{mod} q)$.

Thus for $s \neq 1$ and for $\chi \neq \chi_0 (\operatorname{mod} q)$ , we have

$$L^2(s,\chi) = \sum_{\ell}' \chi(\ell) \sum' \sum_{ab\equiv \ell(q)}' q^{-s} \zeta(s, \tfrac{a}{q}) \cdot q^{-s} \zeta(s, \tfrac{b}{q})$$

$$= \sum_{\ell}' \chi(\ell) \sum' \sum_{ab\equiv \ell(q)}' q^{-2s} \left(\tfrac{q^s}{a^s} + \phi\left(s, \tfrac{a}{q}\right)\right)\left(\tfrac{q^s}{b^s} + \phi\left(s, \tfrac{b}{q}\right)\right) = \sum_{\ell}' \chi(\ell) A(s,\ell,q) \text{, say .}$$

This gives $\sum_{\chi \neq \chi_{0(\operatorname{mod} q)}} |L(s,\chi)|^4 = \sum_{\ell_1}' \sum_{\ell_2}' A(s,\ell_1,q) \overline{A(s,\ell_2,q)} \sum_{\chi \neq \chi_{0(q)}} \chi(\ell_1) \overline{\chi}(\ell_2)$

$$= \phi(q) \sum_{\ell}' |A(s,\ell,q)|^2 - \left| \sum_{\ell}' A(s,\ell,q) \right|^2 \text{ for } s \neq 1 \text{ .}$$

Next , we evaluate $\sum_{\ell}' A(s,\ell,q)$.

We have $\sum_{\ell}' A(s,\ell,q) = \sum_{\ell}' \sum_{ab}' \sum_{\equiv \ell(q)}' q^{-2s} \left(\frac{q^s}{a^s} + \phi\left(s, \frac{a}{q}\right)\right)\left(\frac{q^s}{b^s} + \phi\left(s, \frac{b}{q}\right)\right)$

$$= q^{-2s} \sum_{a}' \left(\tfrac{q^s}{a^s} + \phi\left(s, \tfrac{a}{q}\right)\right) \left( \sum_{\ell}' \sum_{\substack{b \\ ab\equiv \ell(q)}}' \left(\tfrac{q^s}{b^s} + \phi\left(s, \tfrac{b}{q}\right)\right) \right) .$$

Note that for a fixed a , we have $\sum_{\ell}' \sum_{\substack{b \\ ab\equiv \ell(q)}}' \left(\frac{q^s}{b^s} + \phi\left(s, \frac{b}{q}\right)\right) = \sum_{k}' \left(\frac{q^s}{k^s} + \phi\left(s, \frac{k}{q}\right)\right)$ .

Thus we have $\sum_{\ell}' A(s,\ell,q) = \left( \sum_{a}' q^{-s} \left(\frac{q^s}{a^s} + \phi\left(s, \frac{a}{q}\right)\right) \right)^2 = \left( q^{-s} \sum_{a}' \left( \zeta\left(s, \frac{a}{q}\right) - \frac{1}{s-1} \right) \right)^2$

$$= \left( L(s,\chi_0) - \tfrac{\phi(q)}{q^s (s-1)} \right)^2 \text{ for } s \neq 1 \text{ .}$$

This gives $\left| \sum_{\ell}' A(s,\ell,q) \right|^2 = \left| L(s,\chi_0) - \frac{\phi(q)}{q^s (s-1)} \right|^4$ .

Note that $\sum_{\chi \neq \chi_{0(\mathrm{mod}\, q)}} |L(1,\chi)|^4 = \phi(q) \sum_{\ell}{}' |A(1,\ell,a)|^2 - \lim_{s\to 1} \left| \sum_{\ell}{}' A(s,\ell,q) \right|^2$ .

Next $\lim_{s\to 1} \left( \sum_{\ell}{}' A(s,\ell,q) \right)^2 = \left( \lim_{s\to 1} \left( L(s,\chi_0) - \frac{\phi(q)}{q^s(s-1)} \right) \right)^4$ .

Thus we need $\lim_{s\to 1} \left( L(s,\chi_0) - q^{-s} \cdot \frac{\phi(q)}{s-1} \right)$ .

It is enough , if we let $s \to 1$ through real values of $s$ .

Note that $\lim_{s\to 1} \left( L(s,\chi_0) - \frac{\phi(q)}{q^s} \cdot \frac{1}{s-1} \right) = \lim_{s\to 1} \left( \zeta(s) \prod_{p/q} (1-p^{-s}) - \frac{\phi(q)}{q^s(s-1)} \right)$

$$= \lim_{s\to 1} \left( \zeta(s) - \tfrac{1}{s-1} \right) \cdot \prod_{p/q} (1-p^{-s}) + \lim_{s\to 1} \tfrac{1}{s-1} \left( \prod_{p/q} (1-p^{-s}) - \tfrac{\phi(q)}{q^s} \right)$$

$$= \gamma \prod_{p/q} (1-p^{-1}) + \lim_{s\to 1} \tfrac{1}{s-1} \left( \prod_{p/q} (1-p^{-s}) - \tfrac{\phi(q)}{q^s} \right) = \gamma \tfrac{\phi(q)}{q} + \lim_{s\to 1} \tfrac{1}{s-1} \left( \prod_{p/q} (1-p^{-s}) - \tfrac{\phi(q)}{q^s} \right)$$

Next , we obtain $\lim_{s\to 1} \frac{\prod_{p/q} (1-p^{-s}) - \frac{\phi(q)}{q^s}}{s-1}$ . Note that as $s \to 1$ through real values , this limit is of the form $\frac{\prod_{p/q} (1-p^{-1}) - \frac{\phi(q)}{q}}{1-1} = \frac{0}{0}$ .

Thus , using L'Hospital's rule , we have this limit $= \lim_{s\to 1} \frac{d}{ds} \left( \prod_{p/q} (1-p^{-s}) - \frac{\phi(q)}{q^s} \right)$

$$= \lim_{s\to 1} \left\{ \left( \prod_{p/q} (1-p^{-s}) \right) \sum_{p/q} \tfrac{\log p}{p^s - 1} + \tfrac{\phi(q)}{q^s} \log q \right\} = \left( \prod_{p/q} \left(1 - \tfrac{1}{p}\right) \right) \left( \sum_{p/q} \tfrac{\log p}{p-1} \right) + \tfrac{\phi(q)}{q} \log q = \tfrac{\phi(q)}{q} \left( \log q + \sum_{p/q} \tfrac{\log p}{p-1} \right)$$

This gives $\lim_{s\to 1} \left( L(s,\chi_0) - \frac{\phi(q)}{q^s(s-1)} \right) = \frac{\phi(q)}{q} \left( \log q + \sum_{p/q} \frac{\log p}{p-1} + \gamma \right)$.

Thus $\lim_{s\to 1} \left( \sum_{\ell}{}' A(s,\ell,q) \right)^2 = \frac{\phi^4(q)}{q^4} \left( \log q + \sum_{p/q} \frac{\log p}{p-1} + \gamma \right)^4 << \left( \log^4 (q+2) \right) \left( \log\log(q+2) \right)^4$.

Next we estimate $\sum_{\ell}{}' \left( A(1,\ell,q)\right)^2 \; = \sum_{\ell}{}' \left( \sum_{ab}{}' \sum_{\equiv \ell(q)}{}' q^{-2}\left(\tfrac{q}{a}+\phi\left(1,\tfrac{a}{q}\right)\right)\left(\tfrac{q}{b}+\phi\left(1,\tfrac{b}{q}\right)\right)\right)$.

Consider $\sum_{ab}{}' \sum_{\equiv \ell(q)}{}' q^{-2}\left(\tfrac{q}{a}+\phi\left(1,\tfrac{a}{q}\right)\right)\cdot\left(\tfrac{q}{b}+\phi\left(1,\tfrac{b}{q}\right)\right) = \sum_{ab}{}' \sum_{\equiv \ell(q)}{}' \left(\tfrac{1}{a}+\tfrac{1}{q}\phi\left(1,\tfrac{a}{q}\right)\right)\cdot\left(\tfrac{1}{b}+\tfrac{1}{q}\phi\left(1,\tfrac{b}{q}\right)\right)$

$$= \sum{}' \sum_{ab\equiv \ell\ (q)}{}' \left(\tfrac{1}{a}+0(\tfrac{1}{q})\right)\left(\tfrac{1}{b}+0(\tfrac{1}{q})\right)$$

$$= \sum_{ab}{}' \sum_{\equiv \ell(q)}{}' \tfrac{1}{ab} + 0\left(\tfrac{1}{q}\sum_{ab}{}' \sum_{\equiv \ell(q)}{}' \tfrac{1}{a}\right) + 0\left(\tfrac{1}{q}\sum_{ab}{}' \sum_{\equiv \ell(q)}{}' \tfrac{1}{b}\right) + 0\left(\tfrac{1}{q^2}\sum_{a}{}' \sum_{\substack{b \\ ab\equiv \ell(q)}}{}' 1\right)$$

$$= \sum_{ab}{}' \sum_{\equiv \ell(q)}{}' \tfrac{1}{ab} + 0\left(\tfrac{1}{q}\sum_{a}{}' \tfrac{1}{a}\right) + 0\left(\tfrac{1}{q}\sum_{b}{}' \tfrac{1}{b}\right) + 0\left(\tfrac{\phi(q)}{q^2}\right) = \sum_{\substack{n\le q^2 \\ n\equiv \ell(q)}} \tfrac{a(n)}{n} + 0\left(\tfrac{\log(q+3)}{q}\right) .$$

Thus $\phi(q)\sum_{\ell}{}' \left( \sum_{\substack{n\le q^2 \\ n\equiv \ell (\bmod q)}} \tfrac{a(n)}{n} + 0\left(\tfrac{\log(q+3)}{q}\right)\right)^2$

$$= \phi(q)\sum_{\ell}{}' \left( \sum_{\substack{n\le q^2 \\ n\equiv \ell (\bmod q)}} \tfrac{a(n)}{n}\right)^2 + 0\left(\tfrac{\phi(q)}{q}\log(q+3)\cdot \sum_{\substack{n\le q^2 \\ (n,q)=1}} \tfrac{a(n)}{n}\right) + 0\left(\log^2(q+3)\right)$$

$$= \phi(q)\sum_{\ell}{}' \left( \sum_{\substack{n\le q^2 \\ n\equiv \ell (\bmod q)}} \tfrac{a(n)}{n}\right)^2 + 0\left(\log^3(q+3)\right) ,\quad \text{as} \sum_{n\le q^2} \tfrac{a(n)}{n} << \sum_{n\le q^2} \tfrac{d(n)}{n} << \log^2(q+3) .$$

Next , we have $\sum_{\ell}{}' \left( \sum_{\substack{n\le q^2 \\ n\equiv \ell (\bmod q)}} \tfrac{a(n)}{n}\right)^2$

$$= \sum_{\ell}{}' \left( \tfrac{d(\ell)}{\ell} + \sum_{\substack{q<n\le q^2 \\ n\equiv \ell(q)}} \tfrac{a(n)}{n}\right)^2 = \sum_{\ell}{}' \tfrac{d^2(\ell)}{\ell^2} + 2\sum_{\ell}{}' \tfrac{d(\ell)}{\ell}\cdot \sum_{\substack{q<n\le q^2 \\ n\equiv \ell(q)}} \tfrac{a(n)}{n} + \sum_{\ell}{}' \left( \sum_{\substack{q<n\le q^2 \\ n\equiv \ell(q)}} \tfrac{a(n)}{n}\right)^2$$

Next , we estimate $\sum_{\ell}{}^{'} \left( \sum_{\substack{q<n\le q^2 \\ n\equiv \ell (\operatorname{mod} q)}} \frac{a(n)}{n} \right)^2$ .

On using Schwartz's inequality ,

we have $\left( \sum_{\substack{q<n\le q^2 \\ n\equiv \ell (\operatorname{mod} q)}} \frac{a(n)}{n} \right)^2 \le \sum_{\substack{q<n\le q^2 \\ n\equiv \ell (q)}} \frac{1}{n} \cdot \sum_{\substack{q<n\le q^2 \\ n\equiv \ell (q)}} \frac{a^2(n)}{n} \le \sum_{1\le m\le q} \frac{1}{mq+\ell} \cdot \sum_{\substack{q<n\le q^2 \\ n\equiv \ell (q)}} \frac{a^2(n)}{n}$

$\le \sum_{1\le m\le q} \frac{1}{mq} \sum_{\substack{q<n\le q^2 \\ n\equiv \ell (q)}} \frac{a^2(n)}{n} \le \left( \frac{\log(q+3)}{q} \right) \left( \sum_{\substack{q<n\le q^2 \\ n\equiv \ell (q)}} \frac{a^2(n)}{n} \right)$

Thus $\sum_{\ell}{}^{'} \left( \sum_{\substack{q<n\le q^2 \\ n\equiv \ell (q)}} \frac{a(n)}{n} \right)^2 << \left( \frac{\log(q+2)}{q} \right) \sum_{\ell}{}^{'} \sum_{\substack{q<n\le q^2 \\ n\equiv \ell (q)}} \frac{d^2(n)}{n} << \left( \frac{\log(q+2)}{q} \right) \cdot \sum_{q<n\le q^2} \frac{d^2(n)}{n} << \frac{\log^5(q+2)}{q}$

Note $\sum_{\ell}{}^{'} \frac{d^2(\ell)}{\ell^2} << \log^3(q+2)$ .

Hence by Schwartz's lemma , $\sum_{\ell}{}^{'} \frac{d(\ell)}{\ell} \cdot \sum_{\substack{q<n\le q^2 \\ n\equiv \ell (q)}} \frac{a(n)}{n} << \left( \left( \log^3(q+2) \right) \left( \frac{\log^5(q+2)}{q} \right) \right)^{\frac{1}{2}} << \frac{\log^4(q+2)}{\sqrt{q}}$ .

Thus $\phi(q) \sum_{\ell}{}^{'} \left( \sum_{\substack{n\le q^2 \\ n\equiv \ell (q)}} \frac{a(n)}{n} \right)^2 = \phi(q) \sum_{\ell}{}^{'} \frac{d^2(\ell)}{\ell^2} + 0\left( \frac{\phi(q)}{\sqrt{q}} \cdot \log^4(q+2) \right)$ .

Hence $\sum_{\chi \ne \chi_{0(q)}} \left| L(1,\chi) \right|^4 = \phi(q) \sum_{\ell}{}^{'} \frac{d^2(\ell)}{\ell^2} + 0\left( \frac{\phi(q)}{\sqrt{q}} \cdot \log^4(q+2) \right)$ .

This completes the proof of Theorem 1 .

**<u>Proof</u> <u>of</u> <u>Theorem</u> 2** : As in the proof of Theorem 1 , since

$L(s,\chi) = \sum_{a}{}^{'} \chi(a) \left( q^{-s} \cdot \zeta(s, \tfrac{a}{q}) \right)$ , we have for $s \ne 1$ ,

$\sum_{\chi \ne \chi_o (\operatorname{mod} q)} \left| L(s,\chi) \right|^2 = \phi(q) \sum_{a}{}^{'} \left| q^{-s} \cdot \zeta(s, \tfrac{a}{q}) \right|^2 - \left| \sum_{a}{}^{'} q^{-s} \cdot \zeta(s, \tfrac{a}{q}) \right|^2$ .

As in the proof of Theorem 1 , we have $\sum_{\chi\neq\chi_0 (\mathrm{mod}\, q)} |L(1,\chi)|^2$

$$= \phi(q) \sum_a{}' \left| \lim_{s\to 1} q^{-s}\left(\zeta(s,\tfrac{a}{q}) - \tfrac{1}{s-1}\right) \right|^2 - \left| \lim_{s\to 1} \sum_a{}' q^{-s}\left(\zeta(s,\tfrac{a}{q}) - \tfrac{1}{s-1}\right) \right|^2$$

$$= \phi(q) \cdot \sum_a{}' \left| \lim_{s\to 1} q^{-s}\left(\zeta(s,\tfrac{a}{q}) - \tfrac{1}{s-1}\right) \right|^2 - \frac{\phi^2(q)}{q^2}\left( \log q + \sum_{p/q} \frac{\log p}{p-1} + \gamma \right)^2$$

Next , we estimate $\sum_a{}' \left| \lim_{s\to 1} q^{-s}\left(\zeta(s,\tfrac{a}{q}) - \tfrac{1}{s-1}\right) \right|^2$ , using the power series (in $\alpha$ ) of $\zeta_1(s,\alpha)$ as $s \to 1$ . We have $\left| \lim_{s\to 1} q^{-s}\left(\zeta(s,\tfrac{a}{q}) - \tfrac{1}{s-1}\right) \right|^2 = \left| \lim_{s\to 1} q^{-s}\left(\tfrac{q^s}{a^s} + \zeta_1(s,\tfrac{a}{q}) - \tfrac{1}{s-1}\right) \right|^2$

$$= \left| \lim_{s\to 1} a^{-s} + \lim_{s\to 1} q^{-s}(\zeta_1(s,\tfrac{a}{q}) - \tfrac{1}{s-1}) \right|^2 = \left| \tfrac{1}{a} + \tfrac{1}{q}\sum_{n\geq 0} a_n (\tfrac{a}{q})^n \right|^2,$$

where $a_0 = \gamma$ and $a_n = (-1)^n \zeta(n+1)$ for $n \geq 1$ .

Next , $\left( \tfrac{1}{a} + \tfrac{1}{q}\sum_{n\geq 0} a_n (\tfrac{a}{q})^n \right)^2$

$$= \frac{1}{a^s} + \frac{2}{aq}\sum_{n\geq 0} a_n (\tfrac{a}{q})^n + \frac{1}{q^2}\left( \sum_{n\geq 0} a_n (\tfrac{a}{q})^n \right)^2 = \frac{1}{a^2} + \frac{2\gamma}{aq} + \left( \frac{2}{q^2}\sum_{n\geq 1} a_n (\tfrac{a}{q})^{n-1} + \frac{1}{q^2}\sum_{n\geq 0} c_n (\tfrac{a}{q})^n \right)$$

where $c_n = \sum_{i+j=n} a_i a_j$ for $n \geq 0$ .

(Incidentally , note that the expression in the round bracket is $0(\frac{1}{q^2})$ , which is evident in view of the fact that with $N = 1$ in the finite Taylor series development of $\lim_{s\to 1} \left(\zeta_1(s,\alpha) - \frac{1}{s-1}\right)$ , we have $\lim_{s\to 1} q^{-s}\left(\zeta(s,\frac{a}{q}) - \frac{1}{s-1}\right) = \frac{1}{q}\left(\frac{q}{a} + \gamma + 0(\frac{a}{q})\right) = \frac{1}{a} + \frac{\gamma}{q} + 0(\frac{a}{q^2})$ so that $\lim_{s\to 1}\left(q^{-s}\left(\zeta(s,\frac{a}{q}) - \frac{1}{s-1}\right)\right)^2 = \frac{1}{a^2} + \frac{2\gamma}{aq} + 0(\frac{1}{q^2})$ ) .

Thus $\phi(q) \sum_a{}' \left( \frac{1}{a} + \sum_{n\geq 0} a_n (\frac{a}{q})^n \right)^2$

$$=\phi(q){\sum_a}' \tfrac{1}{a^2}+2\gamma\cdot\tfrac{\phi(q)}{q}{\sum_a}'\tfrac{1}{a}+2\ \tfrac{\phi(q)}{q^2}\cdot\sum_{n\geq 1}a_n q^{1-n}{\sum_a}' a^{n-1}+\tfrac{\phi(q)}{q^2}\sum_{n\geq 0}c_n q^{-n}\cdot{\sum_a}' a^n$$

$$=\phi(q){\sum_a}' \tfrac{1}{a^2}+2\gamma\cdot\tfrac{\phi(q)}{q}{\sum_a}'\tfrac{1}{a}+\tfrac{2\phi(q)}{q^2}\sum_{n\geq 1}a_n q^{1-n}\sum_{k/q}\mu(\tfrac{q}{k})\cdot\sum_{a=1}^{k}a^{n-1}+\tfrac{\phi(q)}{q^2}\sum_{n\geq 0}c_n q^{-n}\cdot\sum_{k/q}\mu(\tfrac{q}{k})\cdot\sum_{a=1}^{k}a^n$$

$$=\phi(q){\sum_a}' \tfrac{1}{a^2}+2\gamma\cdot\tfrac{\phi(q)}{q}{\sum_a}'\tfrac{1}{a}$$

$$+\tfrac{2\phi(q)}{q^2}\sum_{n\geq 1}\tfrac{a_n}{n}q^{1-n}\cdot\sum_{k/q}\mu(\tfrac{q}{k})\big(B_n(k+1)-B_n(1)\big)$$

$$+\tfrac{\phi(q)}{q^2}\sum_{n\geq 0}\tfrac{c_n}{n+1}\cdot q^{-n}\cdot\sum_{k/q}\mu(\tfrac{q}{k})\big(B_{n+1}(k+1)-B_{n+1}(1)\big)\ ,$$

where $B_n(u)$ is Bernoulli polynomial of degree $n$.

Thus $\displaystyle\sum_{\chi\neq\chi_o\,(\mathrm{mod}\,q)}|L(1,\chi)|^2=\phi(q){\sum_a}'\left(\tfrac{1}{a}+\sum_{n\geq 0}a_n(\tfrac{a}{q})^n\right)^2-\tfrac{\phi^2(q)}{q^2}\left(\log q+\sum_{p/q}\tfrac{\log p}{p-1}+\gamma\right)^2$

$$=\phi(q){\sum_a}' \tfrac{1}{a^2}+\left(2\gamma\cdot\tfrac{\phi(q)}{q}{\sum_a}'\tfrac{1}{a}-\tfrac{\phi^2(q)}{q^2}\left(\log q+\sum_{p/q}\tfrac{\log p}{p-1}+\gamma\right)^2\right)$$

$$+\tfrac{2\phi(q)}{q^2}\sum_{n\geq 1}\tfrac{a_n q^{1-n}}{n}\sum_{k/q}\mu(\tfrac{q}{k})\big(B_n(k+1)-B_n(1)\big)$$

$$+\tfrac{\phi(q)}{q^2}\sum_{n\geq 0}\tfrac{c_n q^{-n}}{n+1}\cdot\sum_{k/q}\mu(\tfrac{q}{k})\big(B_{n+1}(k+1)-B_{n+1}(1)\big)$$

Next , $\displaystyle 2\gamma\cdot\tfrac{\phi(q)}{q}{\sum_{a=1}^{q}}'\tfrac{1}{a}-\tfrac{\phi^2(q)}{q^2}\left(\log q+\sum_{p/q}\tfrac{\log p}{p-1}+\gamma\right)^2=2\gamma\tfrac{\phi^2(q)}{q^2}\left(\log q+\sum_{p/q}\tfrac{\log p}{p-1}+\gamma\right)+0\left(\tfrac{d(q)\phi(q)}{q^2}\right)$

$$-\tfrac{\phi^2(q)}{q^2}\left(\log q+\sum_{p/q}\tfrac{\log p}{p-1}+\gamma\right)^2$$

$$=-\tfrac{\phi^2(q)}{q^2}\left(\log q+\sum_{p/q}\tfrac{\log p}{p-1}+\gamma\right)\left(\log q+\sum_{p/q}\tfrac{\log p}{p-1}-\gamma\right)+0\left(\tfrac{d(q)\phi(q)}{q^2}\right)$$

$$=-\tfrac{\phi^2(q)}{q^2}\left(\log q+\sum_{p/q}\tfrac{\log p}{p-1}\right)^2+\tfrac{\gamma^2\phi^2(q)}{q^2}+0\left(\tfrac{d(q)\phi(q)}{q^2}\right)$$

$$= \frac{\phi^2(q)}{q^2}\left( \gamma^2 - (\log q + \sum_{p/q} \frac{\log p}{p-1})^2 \right) + 0\left(\frac{d(q)\phi(q)}{q^2}\right)$$

Thus $\displaystyle\sum_{\chi \neq \chi_o \,(\mathrm{mod}\, q)} \left|L(1,\chi)\right|^2 = \phi(q){\sum_a}' \frac{1}{a^2} + \frac{\phi^2(q)}{q^2}\left( \gamma^2 - (\log q + \sum_{p/q} \frac{\log p}{p-1})^2 \right) + 0\left(\frac{d(q)\phi(q)}{q^2}\right)$

$$+ 2\frac{\phi(q)}{q^2}\sum_{n\geq 1} \frac{a_n}{n} q^{1-n} \sum_{k/q} \mu(\tfrac{q}{k})\left(B_n(k+1) - B_n(1)\right)$$

$$+ \frac{\phi(q)}{q^2}\sum_{n\geq o} \frac{c_n q^{-n}}{n+1} \sum_{k/q} \mu(\tfrac{q}{k})\left(B_{n+1}(k+1) - B_{n+1}(1)\right)$$

Incidentally , note that the term $0\left(\frac{d(q)\phi(q)}{q^2}\right)$ can be expressed as a finite asymptotic series with error term $0(q^{-N})$ for any integer $N \geq 1$ .

This completes the proof of Theorem 2 .